\documentclass[]{article}

\usepackage{amsmath} 
\usepackage{graphics} 
\usepackage{epsfig} 
\def\ds{\displaystyle }
\newtheorem{proposition}{Proposition}

\title{Minimization of an energy error functional 
to solve a Cauchy problem arising in plasma physics: 
the reconstruction of the magnetic flux in the vacuum surrounding the plasma in a Tokamak}

\author{Blaise Faugeras$^*$, Amel Ben Abda$^{**}$, \\
	Jacques Blum$^{*}$ and Cedric Boulbe$^{*}$\\
	{\footnotesize $^*$ Laboratoire J.A. Dieudonn\'e, UMR 7351, Universit\'e Nice Sophia Antipolis,}\\
	{\footnotesize 06108 Nice Cedex 02, France}\\
	{\footnotesize $^{**}$ ENIT-LAMSIN, BP 37, 1002 Tunis, Tunisie}\\
	{\footnotesize Email: Blaise.Faugeras@unice.fr}}

\date{}

\begin{document}
\maketitle

\section*{Abstract}
A numerical method for the computation of the magnetic flux 
in the vacuum surrounding the plasma in a Tokamak 
is investigated. It is based on the formulation of a Cauchy problem 
which is solved through the minimization of an energy error functional. 
Several numerical experiments are conducted which show the efficiency of the method.

\section{Introduction}
In order to be able to control the plasma during a fusion experiment in a Tokamak 
it is mandatory to know its position in the vacuum vessel. 
This latter is deduced from the knowledge of the poloidal flux which itself  
relies on measurements of the magnetic field. 
In this paper we investigate a numerical method for the computation of the poloidal flux in the vacuum. 
Let us first briefly recall the equations modelizing the equilibrium of a plasma in a Tokamak 
\cite{Wesson:2004}.

\begin{figure}
\begin{center}
\includegraphics[height=5cm,angle=0]{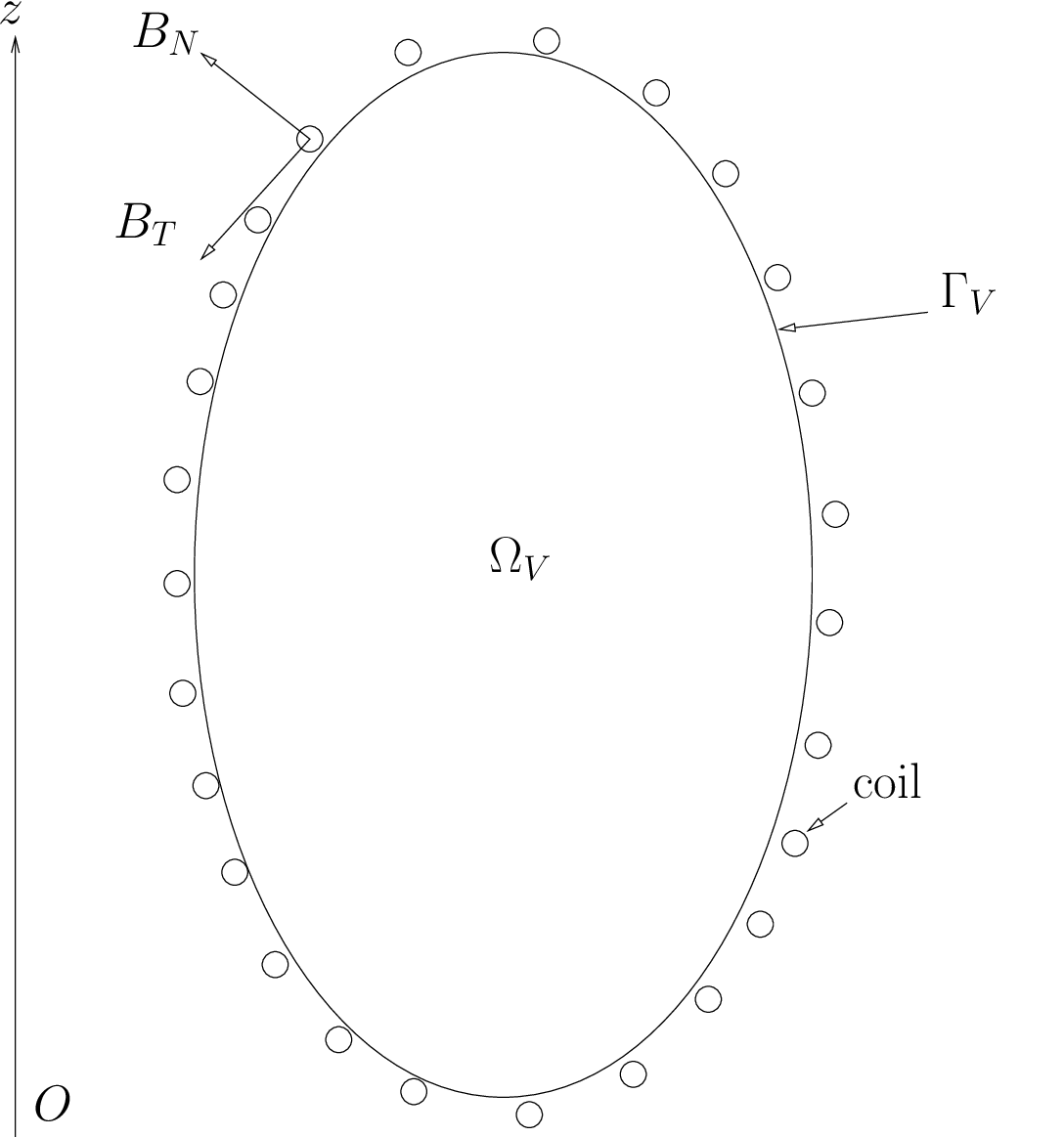}
\caption{\label{fig:fig1} Cross section of the vacuum vessel: 
the domain $\Omega_V$, its boundary $\Gamma_V$. 
Coils providing measurements of the components of the magnetic field tangent and normal 
to $\Gamma_V$ are represented 
surrounding the vacuum vessel.}
\end{center}
\end{figure}

\begin{figure}
\begin{center}
\includegraphics[height=5cm,angle=0]{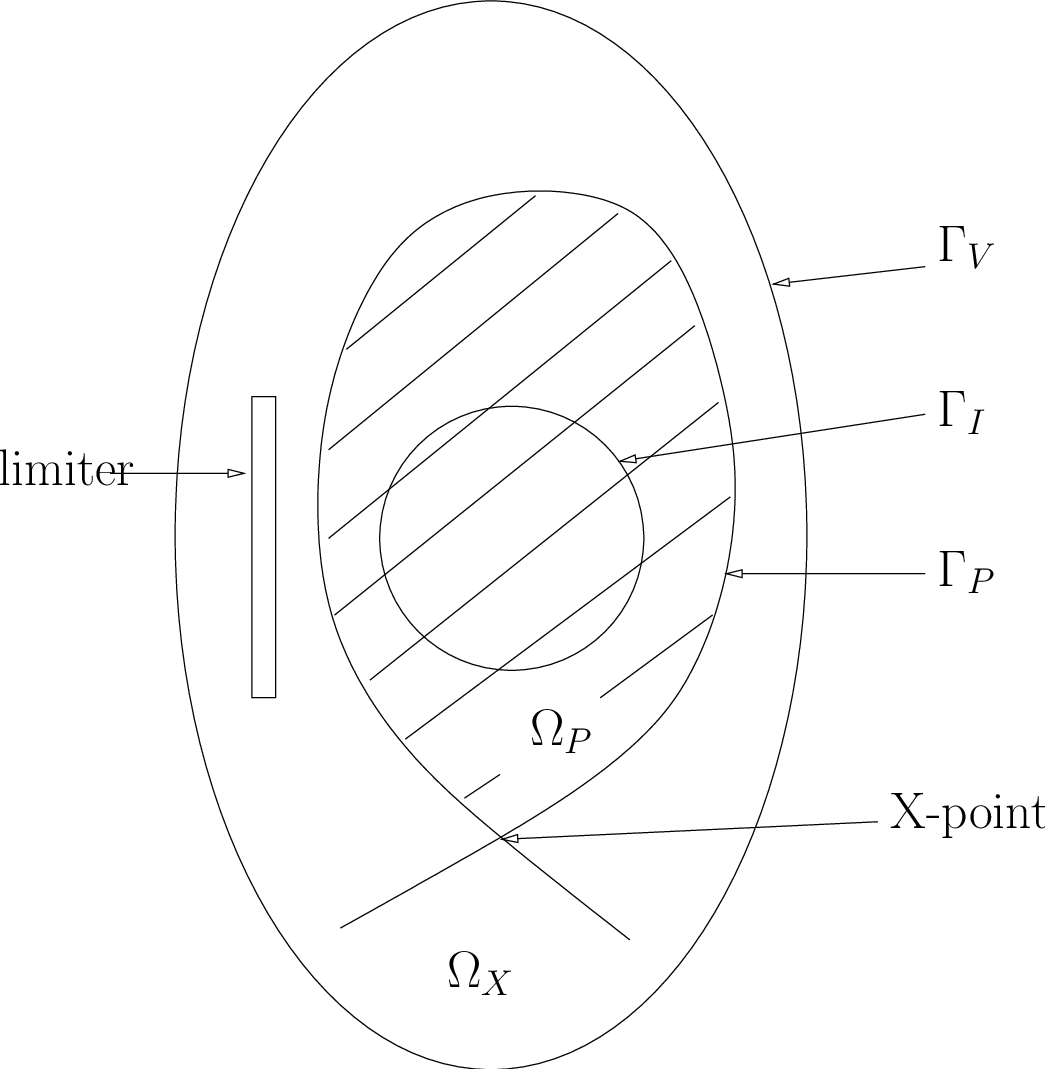}
\includegraphics[height=5cm,angle=0]{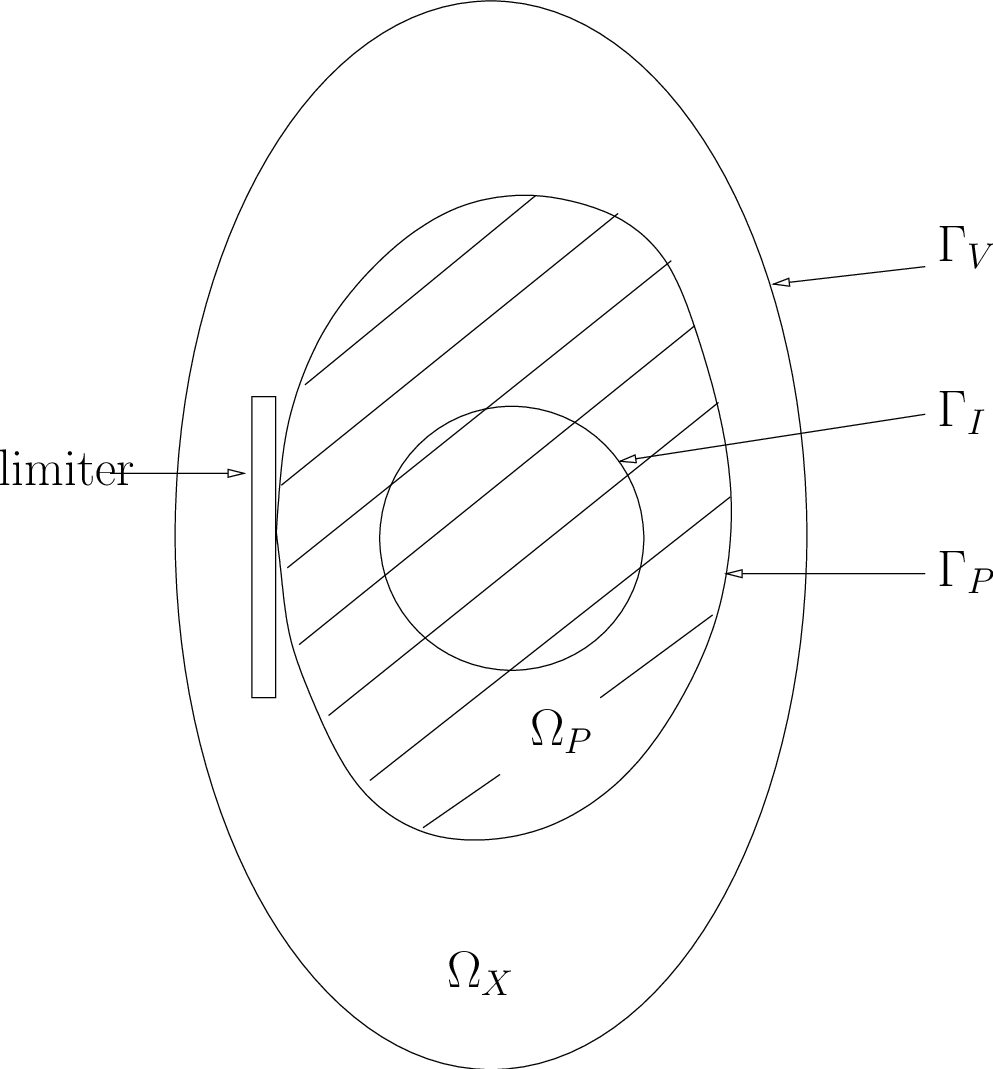}
\caption{\label{fig:fig23} The plasma domain $\Omega_P$ and the vacuum region $\Omega_X$. 
The plasma boundary is determined by an X-point configuration (left) 
or a limiter configuration (right). 
The fictitious contour $\Gamma_I$ is represented inside the plasma.}
\end{center}
\end{figure}

Assuming an axisymmetric configuration one considers 
a 2D poloidal cross section of the 
vacuum vessel $\Omega_V$ in the $(r,z)$ system of coordinates (Fig. \ref{fig:fig1}). 
In this setting the poloidal flux 
$\psi(r,z)$ is related to the magnetic field through 
the relation $(B_r,B_z)=\ds\frac{1}{r}(-\frac{\partial \psi}{\partial z},\frac{\partial \psi}{\partial r})$ 
and, as there is no toroidal current density in the vacuum outside the plasma, satisfies the following equation   
\begin{equation}
\label{eqn:psi}
L \psi=0 \ \mathrm{in}\ \Omega_X
\end{equation}
 where $L$ denotes the elliptic operator 
$$
L.= -[\ds \frac{\partial }{\partial r}(\frac{1}{r}\frac{\partial .}{\partial r}) 
+ \ds \frac{\partial }{\partial z}(\frac{1}{r}\frac{\partial .}{\partial z}) ]
$$ 
and 
$$
\Omega_X=\Omega_V-\bar{\Omega}_P
$$ 
denotes the vacuum region surrounding 
the domain of the plasma $\Omega_P$ of boundary $\Gamma_P$ (see Fig. \ref{fig:fig23}). 
Inside the plasma Eq. (\ref{eqn:psi}) is not valid anymore 
and the poloidal flux satisfies the Grad-Shafranov equation \cite{Shafranov:1958,Grad:1958} 
which describes the equilibrium of a plasma confined by a magnetic field
\begin{equation}
\label{eqn:gs}
L \psi=\mu_0j(r,\psi)\ \mathrm{in}\ \Omega_P
\end{equation}
where $\mu_0$ is the magnetic permeability of the vacuum and $j(r,\psi)$ is the unknown toroidal 
current density function inside the plasma.
Since the plasma boundary $\Gamma_P$ is unknown the equilibrium 
of a plasma in a Tokamak is a free boundary problem described by a particular 
non-linearity of the model. The boundary is an iso-flux line determined 
either as being a magnetic separatrix 
(hyperbolic line with an X-point as on the left hand side of 
Fig. \ref{fig:fig23}) 
or by the contact with a limiter (Fig. \ref{fig:fig23} right hand side). 
In other words the plasma boundary is determined from the equation $\psi(r,z)=\psi_P$, 
$\psi_P$ being the value of the flux at the X-point or the value of the flux for 
the outermost flux line inside a limiter.

In order to compute an approximation of $\psi$ in the vacuum and to find the plasma boundary 
without knowing the current density $j$ in the plasma and thus without using 
the Grad-Shafranov equation (\ref{eqn:gs}) the strategy which is routinely used 
in operational codes mainly consists in choosing an a priori expansion method for $\psi$
such as for example truncated Taylor and Fourier expansions for the code Apolo on 
the Tokamak ToreSupra \cite{Saint-Laurent:2001} 
or piecewise polynomial expansions for the code Xloc on the Tokamak JET 
\cite{OBrien:1993,Sartori:2003}. 
The flux $\psi$ can also be expanded in toroidal harmonics involving Legendre functions 
or expressed by using Green functions in the filament method 
(\cite{Lackner:1976,Feneberg:1984}, \cite{Braams:1991} and the references therein). 
In all cases the coefficients of the
expansion are then computed through a fit to the measurements of the magnetic field. 
Indeed several magnetic probes and flux loops surround the boundary $\Gamma_V$ of the vacuum vessel 
and measure the magnetic field and the flux (see Fig. \ref{fig:fig1}). It should also be noted that very 
similar problems are studied in \cite{Haddar:2005,Bourgeois:2010,Fischer:2011,Fischer:2012}

In this paper we investigate a numerical method based on the resolution 
of a Cauchy problem introduced in (\cite{Blum:1989}, Chapter 5) 
which we recall here below. The proposed approach uses the fact that 
after a preprocessing of these measurements (interpolation and possibly 
integration on a contour) 
one can have access to a complete set of Cauchy data, 
$f=\psi$ on $\Gamma_V$ 
and $g = \ds \frac{1}{r}\frac{\partial \psi}{\partial n}$ on $\Gamma_V$.

The poloidal flux satisfies 
\begin{equation}
\label{eqn:cauchypbm0}
\left \lbrace
\begin{array}{l}
L \psi = 0 \quad \mathrm{in}\ \Omega_X \\[10pt]
\psi= f \quad  \mathrm{on}\  \Gamma_V \\[10pt]
\ds \frac{1}{r}\frac{\partial \psi}{\partial n}= g\quad  \mathrm{on}\  \Gamma_V \\[10pt]
\psi=\psi_P \quad \mathrm{on}\ \Gamma_P 
\end{array} 
\right.
\end{equation}

In this formulation the domain $\Omega_X=\Omega_X(\psi)$ is unknown since 
the free plasma boundary $\Gamma_P$ as well as the flux $\psi_P$ 
on the boundary are unknown. 
Moreover the problem is ill-posed in the sense of Hadamard \cite{Courant:1962} 
since there are two Cauchy conditions on the boundary $\Gamma_V$.
 
In order to compute the flux in the vacuum and to find the plasma boundary we are going 
to define a new problem as in \cite{Blum:1989} which is an approximation of the original one. 
Let us define a fictitious boundary $\Gamma_I$ fixed inside the plasma 
(see Fig. \ref{fig:fig23}). We are going to seek an approximation of the poloidal flux 
$\psi$ satisfying $L\psi=0$ in the domain contained between the fixed boundaries 
$\Gamma_V$ and $\Gamma_I$. The problem becomes one formulated on a fixed domain $\Omega$:
\begin{equation}
\label{eqn:cauchypbm}
\left \lbrace
\begin{array}{l}
L \psi = 0 \quad \mathrm{in}\ \Omega \\[10pt]
\psi= f \quad  \mathrm{on}\  \Gamma_V \\[10pt]
\ds \frac{1}{r}\frac{\partial \psi}{\partial n}= g\quad  \mathrm{on}\  \Gamma_V 
\end{array} 
\right.
\end{equation}

Let us insist here on the fact that this problem is an approximation to the original one since in 
the domain between $\Gamma_P$ and $\Gamma_I$, $\psi$ should satisfy the Grad-Shafranov equation. 
The relevance of this approximating model is consolidated by the Cauchy-Kowalewska theorem \cite{Courant:1962}. 
For $\Gamma_P$ smooth enough the function $\psi$ can be extended in the sense 
of $L\psi=0$ in a neighborhood of $\Gamma_P$ inside the plasma. 
Hence the problem formulated on a fixed domain with a fictitious boundary 
$\Gamma_I$ not "too far" from $\Gamma_P$ is an approximation of the free boundary problem. 
As mentioned in \cite{Blum:1989} if $\Gamma_I$ were identical with $\Gamma_P$ 
then by the virtual shell principle \cite{Shafranov:1972} 
the quantity $w=\ds \frac{1}{r}\frac{\partial \psi}{\partial n}|_{\Gamma_I}$ would 
represent the surface current density (up to a factor $\ds \frac{1}{\mu_0}$) 
on $\Gamma_P$ for which the magnetic field created outside the plasma by the current 
sheet is identical to the field created by the real current density spread throughout the plasma. 

However no boundary condition is known on $\Gamma_I$. 
One way to deal with this second issue and to solve such a problem 
is to formulate it as an optimal control one. 
Only the Dirichlet condition on $\Gamma_V$ is retained 
to solve the boundary value problem and a least square error functional measuring 
the distance between measured and computed normal derivative and depending 
on the unknown boundary condition on $\Gamma_I$ is minimized. 
Due to the illposedness of the considered Cauchy problem a regularization 
term is needed to avoid erratic behaviour on the boundary where the data is missing. 
A drawback of this method developed in \cite{Blum:1989} is that Dirichlet 
and Neumann boundary conditions on $\Gamma_V$ are not used in a symmetric way. 
One is used as a boundary condition for the partial differential equation, $L\psi=0$, whereas the other is 
used in the functional to be minimized.

Freezing the domain to $\Omega$ by introducing the fictitious boundary $\Gamma_I$ 
enables to remove the nonlinearity of the problem. The plasma boundary $\Gamma_P$ 
can still be computed as an iso-flux line and thus is an output of our computations. 
We are going to compute a function $\psi$ such that 
the Dirichlet boundary condition $u=\psi$ on $\Gamma_I$ is such 
that the Cauchy conditions on $\Gamma_V$ are satisfied 
as nearly as possible in the sense of the error functional defined in the next Section.

The originality of the approach proposed in this paper 
relies on the use of an error functional having a physical meaning: 
an energy error functional or constitutive law error functional. 
Up to our knowledge this misfit functional has been introduced in \cite{Ladeveze:1983} 
in the context of a posteriori estimator in the finite element method. 
In this context, the minimization of the constitutive law error functional 
allows to detect the reliability of the mesh 
without knowing the exact solution.
Within the inverse problem community this functional has been introduced in 
\cite{Kohn:1985, Kohn:1989, Kohn:1990} in the context of parameter identification. 
It has been widely exploited in the same context in \cite{Bonnet:2005}. 
It has also been used for Robin type boundary condition recovering \cite{Chaabane:1999} and 
in the context of geometrical flaws identification (see \cite{BenAbda:2009} and references therein). 
For lacking boundary data recovering (i.e. Cauchy problem resolution) 
in the context of Laplace operator, 
the energy error functional has been introduced in \cite{Andrieux:2005, Andrieux:2006}.  
A study of similar techniques can be found in \cite{BenBelgacem:2005, Azaiez:2006} and the analysis found 
in these papers uses elements taken 
from the domain decomposition framework \cite{Quarteroni:1999}.

The paper is organized as follows. 
In Section 2 we give the formulation of the problem we are interested in 
and provide an analysis of its well posedness. 
Section 3 describes the numerical method used. 
Several numerical experiments are conducted to validate it. 
The final experiment shows the reconstruction of the poloidal 
flux and the localization of the plasma boundary 
for an ITER configuration. 

\section{Formulation and analysis of the method}

\subsection{Problem formulation}

As described in the Introduction the starting point is the free boundary problem 
(\ref{eqn:cauchypbm0}). 
We first proceed as in \cite{Blum:1989} and in a first step consider 
the fictitious contour $\Gamma_I$ fixed in the plasma 
and the fixed domain $\Omega$ contained 
between $\Gamma_V$ and $\Gamma_I$. Problem 
(\ref{eqn:cauchypbm0}) is approximated by the Cauchy problem (\ref{eqn:cauchypbm}). 
The boundaries $\Gamma_V$ and $\Gamma_I$ are assumed to be chosen smooth enough in order not 
to refrain any of the developments which follow in the paper.

In a second step the problem is separated into two different ones. 
In the first one we retain the Dirichlet boundary condition on $\Gamma_V$ only, 
assume $v$ is given 
on $\Gamma_I$ and seek the solution $\psi_D$ 
of the well-posed boundary value problem:

\begin{equation}
\label{eqn:psid}
\left\{
\begin{array}{l}
L \psi_D = 0 \quad \mathrm{in}\ \Omega \\[10pt]
\psi_D = f \quad  \mathrm{on}\  \Gamma_V \\[10pt]
\psi_D = v\quad  \mathrm{on}\  \Gamma_I 
\end{array}
\right.
\end{equation}

The solution $\psi_D$ 
can be decomposed  in a part linearly depending on $v$ 
and a part depending on $f$ only. 
We have the following decomposition:
\begin{equation}
\label{eqn:decompositionD}
\psi_D=\psi_D(v,f)=\psi_D(v,0)+\psi_D(0,f):=\psi_D(v)+\tilde{\psi}_D(f)
\end{equation}

where $\psi_D(v)$ and $\tilde{\psi}_D(f)$ satisfy:
\begin{equation}
\label{eqn:decomposition_Psid}
\begin{array}{lll}
\left\{
\begin{array}{l}
L\psi_D(v) = 0 \quad \mathrm{in}\ \Omega \\[10pt]
\psi_D(v) = 0 \quad  \mathrm{on}\  \Gamma_V \\[10pt]
\psi_D(v) = v\quad  \mathrm{on}\  \Gamma_I 
\end{array}
\right.
&\hskip 1cm  \hskip 1cm &
\left\{
\begin{array}{l}
L\tilde{\psi}_D(f) = 0 \quad \mathrm{in}\ \Omega \\[10pt]
\tilde{\psi}_D(f)= f\quad  \mathrm{on}\  \Gamma_V \\[10pt]
\tilde{\psi}_D(f)= 0 \quad  \mathrm{on}\  \Gamma_I 
\end{array}
\right.
\end{array}
\end{equation}

In the second problem we retain the Neumann boundary condition only 
and look for $\psi_N$ satisfying the well-posed boundary value problem:
\begin{equation}
\label{eqn:psin}
\left\{
\begin{array}{l}
L \psi_N = 0 \quad \mathrm{in}\ \Omega \\[10pt]
\ds \frac{1}{r}\frac{\partial \psi_N}{\partial n}= g\quad  \mathrm{on}\  \Gamma_V \\[10pt]
\psi_N= v \quad  \mathrm{on}\  \Gamma_I 
\end{array}
\right.
\end{equation}

in which $\psi_N$  can be decomposed  in a part linearly depending on $v$ 
and a part depending on $g$ only. 
We have the following decomposition:
\begin{equation}
\label{eqn:decompositionN}
\psi_N=\psi_N(v,g)=\psi_N(v,0)+\psi_N(0,g):=\psi_N(v)+\tilde{\psi}_N(g)
\end{equation}
where
\begin{equation}\label{eqn:decomposition_Psin}
\begin{array}{lll}
\left\{
\begin{array}{l}
L\psi_N(v) = 0 \quad \mathrm{in}\ \Omega \\[10pt]
\ds \frac{1}{r}\frac{\partial\psi_N(v)}{\partial n} = 0 \quad  \mathrm{on}\  \Gamma_V \\[10pt]
\psi_N(v) = v\quad  \mathrm{on}\  \Gamma_I 
\end{array}
\right.
&\hskip 1cm  \hskip 1cm &
\left\{
\begin{array}{l}
L\tilde{\psi}_N(g) = 0 \quad \mathrm{in}\ \Omega \\[10pt]
\ds \frac{1}{r}\frac{\partial \tilde{\psi}_N}{\partial n}= g\quad  \mathrm{on}\  \Gamma_V \\[10pt]
\psi_N= 0 \quad  \mathrm{on}\  \Gamma_I 
\end{array}
\right.
\end{array}
\end{equation}

In order to solve problem (\ref{eqn:cauchypbm}), 
$f \in H^{1/2}(\Gamma_V)$ and $g \in H^{-1/2}(\Gamma_V)$ 
being given, we would like to find 
$u \in \mathcal{U} = H^{1/2}(\Gamma_I)$ such that $\psi=\psi_D(u,f)=\psi_N(u,g)$. 
To achieve this we are in fact going to seek $u$ 
such that $J(u)=\underset{v \in \mathcal{U}}{\mathrm{inf}} J(v)$ 
where $J$ is the error functional defined by
\begin{equation}
\label{eqn:kohnvog}
J(u)=\ds \frac{1}{2} \int_{\Omega} \frac{1}{r} 
|| \nabla \psi_D(u,f) - \nabla \psi_N(u,g) ||^2 dx 
\end{equation}
measuring a misfit between the Dirichlet solution and the Neumann solution.

\subsection{Analysis of the method}

In order to minimize $J$ one can compute its derivative 
and express the first order optimality condition. 
When doing so the two symmetric bilinear forms $s_D$ and $s_N$ 
as well as the linear form $l$ defined below appear naturally 
and in a first step it is convenient to give 
a new expression of functional (\ref{eqn:kohnvog}) 
using these forms.

Let $u,v \in H^{1/2}(\Gamma_I)$ and define
\begin{equation}
s_D(u,v)= \ds \int_\Omega \ds \frac{1}{r} \nabla \psi_D(u) \nabla \psi_D(v) dx 
\end{equation}
Applying Green's formula and noticing that 
$\psi_D(v)=v$ on $\Gamma_I$ and $\psi_D(v)=0$ 
on $\Gamma_V$ we obtain
\begin{equation}
\label{eqn:sdnormal}
s_D(u,v)= 
\ds \int_{\partial \Omega} \ds \frac{1}{r} \partial_n \psi_D(u) \psi_D(v) d\sigma - \ds \int_\Omega \nabla (\ds \frac{1}{r} \nabla \psi_D(u)) \psi_D(v) dx
=\ds \int_{\Gamma_I} \ds \frac{1}{r}  \partial_n \psi_D(u) v d\sigma  
\end{equation}
where the integrals on the boundary are to be understood as duality pairings. 
In Eq. (\ref{eqn:sdnormal}) one can replace $\psi_D(v)$ by any extension 
$\mathcal{R}(v)$ in $H^1_0(\Omega,\Gamma_V)= \lbrace \psi \in H^1(\Omega), \psi_{|\Gamma_V}=0 \rbrace$
of $v \in H^{1/2}(\Gamma_I)$.\\ 
Hence $s_D$ can be represented by
\begin{equation}
s_D(u,v)= \ds \int_\Omega \ds \frac{1}{r} \nabla \psi_D(u) \nabla \mathcal{R}(v) dx 
\end{equation}

Equivalently $s_N$ is defined by
\begin{equation}
s_N(u,v)= \ds \int_\Omega \ds \frac{1}{r} \nabla \psi_N(u) \nabla \psi_N(v) dx 
\end{equation}
Since $\psi_N(v)=v$ on $\Gamma_I$ and $\ds \frac{1}{r}\partial_n \psi_N(u)=0$ 
on $\Gamma_V$ we have that 
\begin{equation}
\label{eqn:snnormal}
s_N(u,v)=
\ds \int_{\partial \Omega} \ds \frac{1}{r} \partial_n \psi_N(u) \psi_N(v) d\sigma - \ds \int_\Omega \nabla (\ds \frac{1}{r} \nabla \psi_N(u)) \psi_N(v) dx
=\ds \int_{\Gamma_I}\ds \frac{1}{r} \partial_n \psi_N(u) v d\sigma 
\end{equation}
and $s_N$ can also be represented by
\begin{equation}
s_N(u,v)= \ds \int_\Omega \ds \frac{1}{r} \nabla \psi_N(u) \nabla \mathcal{R}(v) dx 
\end{equation}
where $\mathcal{R}(v)$ is any extension in $H^1(\Omega)$ of $v \in H^{1/2}(\Gamma_I)$.

Let us now introduce
\begin{equation}
F(u,v)=\ds \frac{1}{2} \int_\Omega \ds \frac{1}{r}(\nabla \psi_D(u,f)-\nabla \psi_N(u,g))
(\nabla \psi_D(v,f)-\nabla \psi_N(v,g))dx
\end{equation}
such that $J(v)=F(v,v)$ 
and the linear form $l$ defined by
\begin{equation}
l(v)=-\ds \int_\Omega \ds \frac{1}{r} (\nabla \tilde{\psi}_D(f) - 
\nabla \tilde{\psi}_N(g))\nabla \psi_D(v)dx
\end{equation}
which can also be computed as
\begin{equation}
l(v)=-\ds \int_\Omega \ds \frac{1}{r} (\nabla \tilde{\psi}_D(f) - 
\nabla \tilde{\psi}_N(g))\nabla \mathcal{R}(v)dx
\end{equation}

It can then be shown that 
\begin{equation}
F(u,v)=\ds \frac{1}{2}(s_D(u,v)-s_N(u,v)-l(u)-l(v))+c 
\end{equation}
where the constant $c$ is given by
\begin{equation}
c=\ds \frac{1}{2} \int_\Omega \ds \frac{1}{r} ||\nabla \tilde{\psi}_D(f) - \nabla \tilde{\psi}_N(g)||^2dx
\end{equation}

Hence functional $J$ can be rewritten as
\begin{equation}
J(v)=\ds \frac{1}{2}(s_D(v,v)- s_N(v,v) )-l(v)+c 
\end{equation}

Following the analysis provided in \cite{BenBelgacem:2005} it can be proved that in the 
favorable case of compatible Cauchy data $(f,g)$ the Cauchy problem admits a solution. 
There exists a unique $u \in \mathcal{U}$ such that $\psi_D(u,f)=\psi_N(u,g)$. 
The minimum of $J$ is also uniquely reached at this point, $J(u)=0$. 
This solution is given by the first order 
optimality condition which reads
\begin{equation}
\label{eqn:opt0}
(J'(u),v)= s_D(u,v) - s_N(u,v) -l(v)=0\quad \forall v \in \mathcal{U}
\end{equation}
Equation (\ref{eqn:opt0}) has an interpretation in terms 
of the normal derivative of $\psi_D$ and $\psi_N$ on the boundary. 
From Eqs. (\ref{eqn:sdnormal}) and (\ref{eqn:snnormal}) and from
\begin{equation}
l(v)=-\ds \int_\Omega \ds \frac{1}{r} (\nabla \tilde{\psi}_D(f) - 
\nabla \tilde{\psi}_N(g))\nabla \psi_D(v)dx
=-\ds \int_{\Gamma_I} \frac{1}{r}(\partial_n \tilde{\psi}_D(f) 
- \partial_n \tilde{\psi}_N(g)) v d\sigma
\end{equation}
we deduce that the optimality condition can be rewritten as
\begin{equation}
\label{eqn:normal-boundary}
\ds \int_{\Gamma_I} 
[(\frac{1}{r}\partial_n \psi_D(u,f) 
- \frac{1}{r}\partial_n \psi_N(u,g))] v d\sigma =0 \quad 
\forall v \in \mathcal{U}
\end{equation}
which can be understood as the equality of the normal derivatives on $\Gamma_I$.

Hence the 
first optimality condition when minimizing $J$ amounts to solve 
an interfacial equation
$$
(S_{D}-S_{N})(v)=\chi,
$$
where $S_{D}$ and $S_{N}$ are the Dirichlet-to-Neumann operators associated to 
the bilinear forms and defined by:
 
\begin{equation}\label{eq31}
\begin{array}{ccccc}
  S_{D} &:& H^{1/2}(\Gamma_I)& \longrightarrow & H^{-1/2}(\Gamma_I)\\
      & & v &\longrightarrow  & \displaystyle\frac{1}{r}
\displaystyle\frac{\partial{\psi}_D(v)}{\partial n}.
\end{array}
\end{equation}
\begin{equation}\label{eq32}
\begin{array}{ccccc}
  S_{N} &:& H^{1/2}(\Gamma_I)& \longrightarrow & H^{-1/2}(\Gamma_I)\\
      & & v &\longrightarrow  & \displaystyle\frac{1}{r}\displaystyle\frac{\partial{\psi}_N(v)}{\partial n},
\end{array}
\end{equation}
and $\chi=-\displaystyle\frac{1}{r}\displaystyle\frac{\partial\tilde{\psi}_D}{\partial n}+\displaystyle\frac{1}{r}
\displaystyle\frac{\partial\tilde{\psi}_N}{\partial n}\ \mbox{on}\ \Gamma_{I}.
$

%
%

Since $S_{D}$ and $S_{N}$ have the same eigenvectors and have asymptotically the same eigenvalues, 
the interfacial operator $S=S_{D}-S_{N}$ is almost singular \cite{BenBelgacem:2005}. 
This point together with the fact that the set of incompatible Cauchy data is known to be dense in 
the set of compatible data (and thus numerical Cauchy data can hardly by compatible) make 
this inverse problem severely ill-posed. 

Some regularization process has to be used. One way to regularize the problem 
is to directly deal with the resolution of the underlying quasi-singular linear system 
using for example a relaxed gradient method \cite{Andrieux:2005, Andrieux:2006}. 
In this paper we have chosen a regularization method of the Tikhonov type. 
It consists in shifting the spectrum of $S$ by adding a term
$$
(S_{D}-S_{N})+\varepsilon S_{D}.
$$ 
where $\varepsilon$ is a small regularization parameter. 
This regularization method is quite natural since the ill-posedness 
of the inverse problem and the lack of stability 
in the identification of $u$ by the minimization of $J$ 
is strongly linked to the fact that $J$ is not coercive 
(see \cite{BenBelgacem:2005} and below).  
We are thus going to minimize the regularized cost function:
$$
J_\varepsilon(v)=J(v)+\varepsilon R_D(v)
$$
with 
$$
R_D(v)=\ds \frac{1}{2} \int_\Omega \frac{1}{r} ||\nabla \psi_D(v)||^2 dx
$$

This brings us to the framework described in \cite{Lions2:1968}. 
We want to solve the following 

\begin{center}
Problem $P_\varepsilon$: $\quad $  find $u_\varepsilon \in \mathcal{U}$ 
such that $J_\varepsilon(u_\varepsilon)
=\underset{v \in \mathcal{U}}{\mathrm{inf}} J_\varepsilon(v)$
\end{center}
 
and the following result holds.

\begin{proposition}

\begin{enumerate}

\item Problem $P_\varepsilon$ admits a unique solution $u_\varepsilon \in \mathcal{U}$ 
characterized by the first order optimality condition 
\begin{equation}
\label{eqn:opt00}
(J'_\varepsilon(u_\varepsilon),v)=\varepsilon s_D(u_\varepsilon,v) + s_D(u_\varepsilon,v) - s_N(u_\varepsilon,v) 
-l(v)=0\quad \forall v \in \mathcal{U}
\end{equation}

\item For a fixed $\varepsilon$ the solution is stable with respect to the data $f$ and $g$. \\
If $f^1$, $f^2 \in H^{1/2}(\Gamma_V)$ and $g^1$, $g^2 \in H^{-1/2}(\Gamma_V)$ it holds that
\begin{equation}
|| u^1_\varepsilon - u^2_\varepsilon ||_{H^{1/2}(\Gamma_I)} \le \ds \frac{C}{\varepsilon} (|| f^1 - f^2 ||_{H^{1/2}(\Gamma_V)} 
+ ||g^1 - g^2 ||_{H^{-1/2}(\Gamma_V)})
\end{equation}

\item If there exists $u \in \mathcal{U}$ such that $\psi_D(u,f)=\psi_N(u,g)$ then 
$u_\varepsilon \rightarrow u$ in $\mathcal{U}$ when $\varepsilon \rightarrow 0$. 

\end{enumerate}

\end{proposition}

Elements of the proof are given in Appendix.

\section{Numerical method and experiments}

\subsection{Finite element discretization}

The resolution of the boundary value 
problems (\ref{eqn:decomposition_Psid}) and (\ref{eqn:decomposition_Psin}) is based on 
a classical $P^1$ finite element method \cite{Ciarlet:1980}.

Let us consider the family of triangulation ${\tau}_h$ of
$\Omega$, and $V_h$ the finite dimensional subspace of $H^1(\Omega)$
defined by
$$
V_h=\{\psi_h\in H^1(\Omega), \psi_{h|T}\in P^1(T),\,\forall T\in {\tau}_h\}.
$$
Let us also introduce the finite element space on $\Gamma_I$
$$
D_h=\{v_h=\psi_h|_{\Gamma_I},\  \psi_h \in V_h\}.
$$
Consider $(\phi_i)_{i=1,...N}$ a basis of $V_h$ and assume that the first $N_{\Gamma_I}$ 
mesh nodes (and basis functions) correspond to the ones situated on $\Gamma_I$. 
A function $\psi_h \in V_h$ is decomposed as $\psi_h=\sum_{i=1}^N a_i \phi_i$ and its trace on 
$\Gamma_I$ as $v_h=\psi_h|_{\Gamma_I}=\sum_{i=1}^{N_{\Gamma_I}} a_i \phi_i|_{\Gamma_I}$.

Given boundary conditions $v_h$ on $\Gamma_I$  and $f_h$, $g_h$ on $\Gamma_V$ 
one can compute the approximations 
$\psi_{D,h}(v_h)$,  $\psi_{N,h}(v_h)$, $\tilde{\psi}_{D,h}(f_h)$ and $\tilde{\psi}_{N,h}(g_h)$ 
with the finite element method.

In order to minimize the discrete regularized error functional, $J_{\varepsilon,h}(u_h)$ 
we have to solve the discrete optimality condition which reads
\begin{equation}
\label{eqn:optdiscrete}
\varepsilon s_{D,h}(u_{h},v_h) + s_{D,h}(u_h,v_h) - s_{N,h}(u_h,v_h) - l(v_h)=0\quad \forall v_h \in D_h
\end{equation}
which is equivalent to look for the vector $\mathbf{u}$ solution to the linear system
\begin{equation}
\label{eqn:linearsystem}
 \mathbf{Su}=\mathbf{l}
\end{equation}
where the $N_{\Gamma_I} \times N_{\Gamma_I}$ matrix $\mathbf{S}$ representing 
the bilinear form $s_h=\varepsilon s_{D,h} + s_{D,h} -s_{N,h}$ 
is defined by
\begin{equation}
\mathbf{S}_{ij}= s_h(\phi_i,\phi_j)
\end{equation}
and $\mathbf{l}$ is the vector $(l_h(\phi_i))_{i=1,...N_{\Gamma_I}}$.

In order to lighten the computations the matrices are evaluated by
\begin{equation}
s_{D,h}(\phi_i,\phi_j)=\ds \int_\Omega \frac{1}{r}\nabla \psi_{D,h}(\phi_i) \nabla \mathcal{R} (\phi_j) dx
\end{equation} 
and 
\begin{equation}
s_{N,h}(\phi_i,\phi_j)=\ds \int_\Omega \frac{1}{r}\nabla \psi_{N,h}(\phi_i) \nabla \mathcal{R} (\phi_j) dx
\end{equation} 
where $\mathcal{R}(\phi_j)$ is the trivial extension which coincides with $\phi_j$ on $\Gamma_I$ 
and vanishes elsewhere.

In the same way the right hand side $\mathbf{l}$ is evaluated by
\begin{equation}
l_h(\phi_i)=-\ds \int_\Omega \frac{1}{r} ( \nabla \tilde{\psi}_{D,h}(f_h)- \nabla \tilde{\psi}_{N,h}(g_h)) 
\nabla \mathcal{R} (\phi_i) dx
\end{equation}

It should be noticed here that matrix $\mathbf{S}$ depends on the geometry of the problem only 
and not on the input Cauchy data. 
Therefore it can be computed once for all  (as well as its $LU$ decomposition for exemple if this is the method used to invert the system) 
and be used for the resolution of successive problems with varying input data 
as it is the case during a plasma shot in a Tokamak. 
Only the right hand side $\mathbf{l}$ has to be recomputed. 
This enables very fast computation times.

\begin{figure}[h]
\begin{center}
\includegraphics[height=6cm,angle=0]{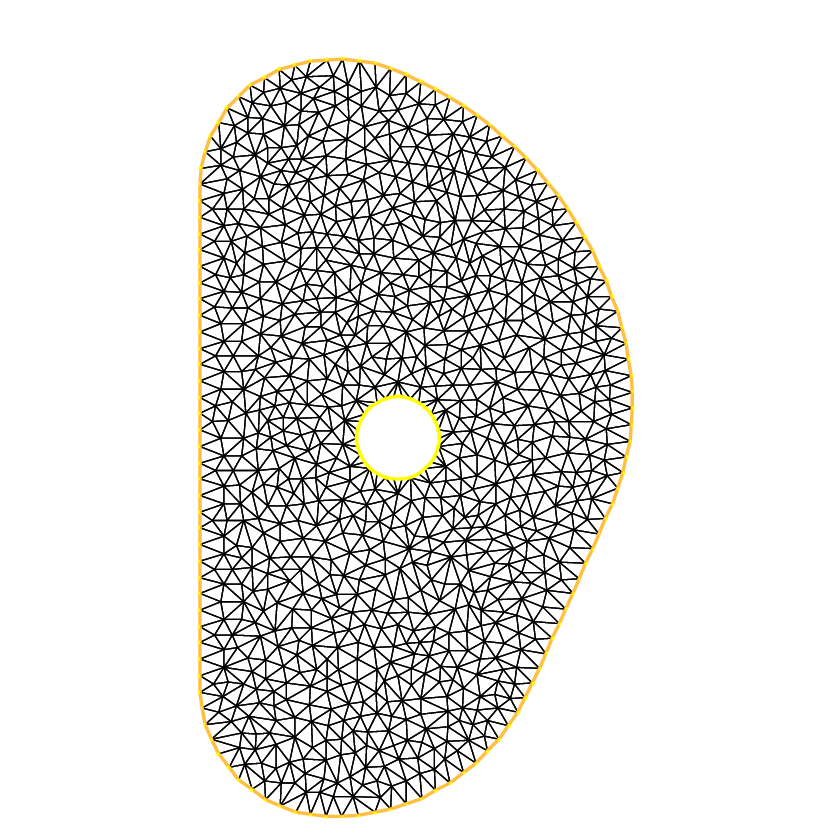}
\caption{\label{fig:itermesh}The mesh used for the ITER configuration in FreeFem++}
\end{center}
\end{figure}

All the numerical results presented in the remaining part of this paper were obtained using the software FreeFem++ (http://www.freefem.org/ff++/). We are concerned with the geometry of ITER and the mesh used for the computations is shown on Fig. \ref{fig:itermesh}. 
It is composed of 1804 triangles and 977 nodes 150 of which are boundary nodes divided into 
120 nodes on $\Gamma_V$ and 30$=N_{\Gamma_I}$ on $\Gamma_I$. 
The shape of $\Gamma_I$ is chosen empirically. 

\subsection{Twin experiments}
Numerical experiments with simulated input Cauchy data are conducted in order to validate the algorithm.
Assume we are provided with a Neumann boundary condition 
function $g$ on $\Gamma_V$. We generate the associated Dirichlet 
function $f$ on $\Gamma_V$ assuming a reference Dirichlet 
function $u_{ref}$ is known on $\Gamma_I$. We thus solve the following boundary value problem:
 
\begin{equation}
\label{eqn:twinref}
\left\{
\begin{array}{l}
L\psi_{N,ref}(u_{ref},g) = 0 \quad \mathrm{in}\ \Omega \\[10pt]
\ds \frac{1}{r} \partial_n\psi_{N,ref}(u_{ref},g) = g \quad  \mathrm{on}\  \Gamma_V \\[10pt]
\psi_{N,ref}(u_{ref},g) = u_{ref} \quad  \mathrm{on}\  \Gamma_I
\end{array}
\right.
\end{equation}

and set $f=\psi_{N,ref}(u_{ref},g)|_{\Gamma_V}$.

We have considered two test cases. In the first one (TC1)
\begin{equation}
\label{eqn:uref1} 
u_{ref}(r,z)=50 \sin(r)^2 + 50 \quad \mathrm{on}\ \Gamma_I
\end{equation}
and in the second one (TC2) $u_{ref}$ is simply a constant
\begin{equation}
\label{eqn:uref2} 
u_{ref}(r,z)=40 \quad \mathrm{on}\ \Gamma_I
\end{equation}

The numerical experiments consist in minimizing the regularized error functional 
$J_\varepsilon$ defined thanks to $f$ and $g$. The obtained optimal solution $u_{opt}$ and 
the associated $\psi_{opt}$ are then compared to $u_{ref}$ and $\psi_{ref}$ which should ideally be recovered. Three cases are considered: the noise free case, 
a $1\%$ noise on $f$ and $g$ and a $5\%$ noise. 

When the noise on $f$ and $g$ is small and the recovery of $u$ is excellent there is 
very little difference between the Dirichlet solution $\psi_D(u_{opt},f)$ and the Neumann solution 
$\psi_N(u_{opt},g)$. However this is not the case any longer 
when the level of noise increases. The Dirichlet solution is much more sensitive to noise 
on $f$ than the Neumann solution is sensitive to noise on $g$. Therefore the optimal solution is chosen to be 
$\psi_{opt}=\psi_N(u_{opt},g)$.

\begin{figure}
\begin{center}
\begin{tabular}{ll}
\includegraphics[height=6cm,angle=0]{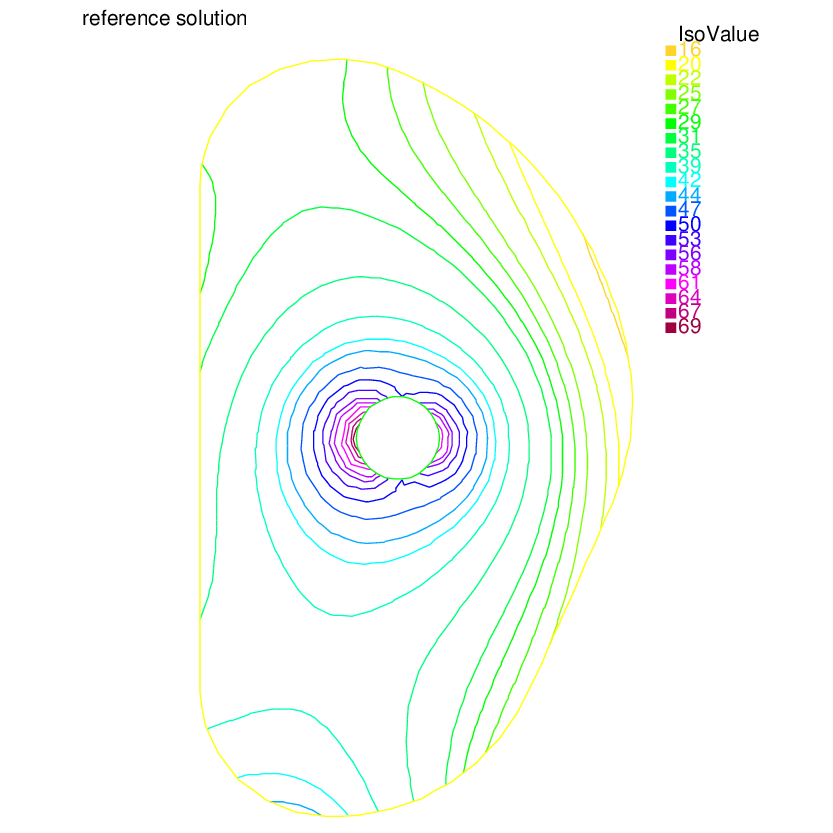}& 
\includegraphics[height=6cm,angle=0]{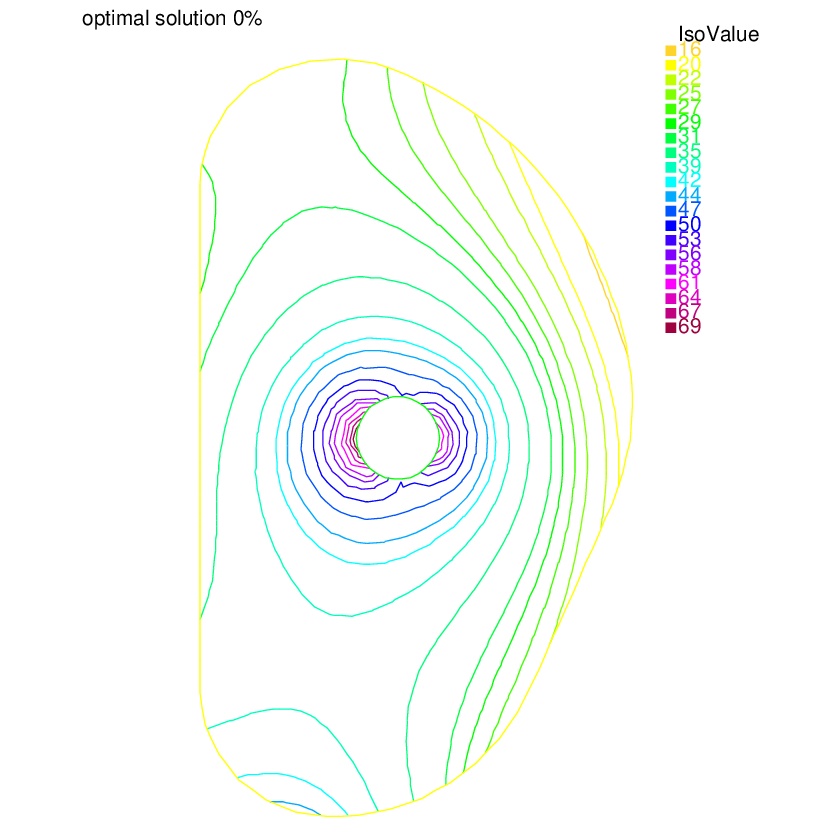}\\ 
\includegraphics[height=6cm,angle=0]{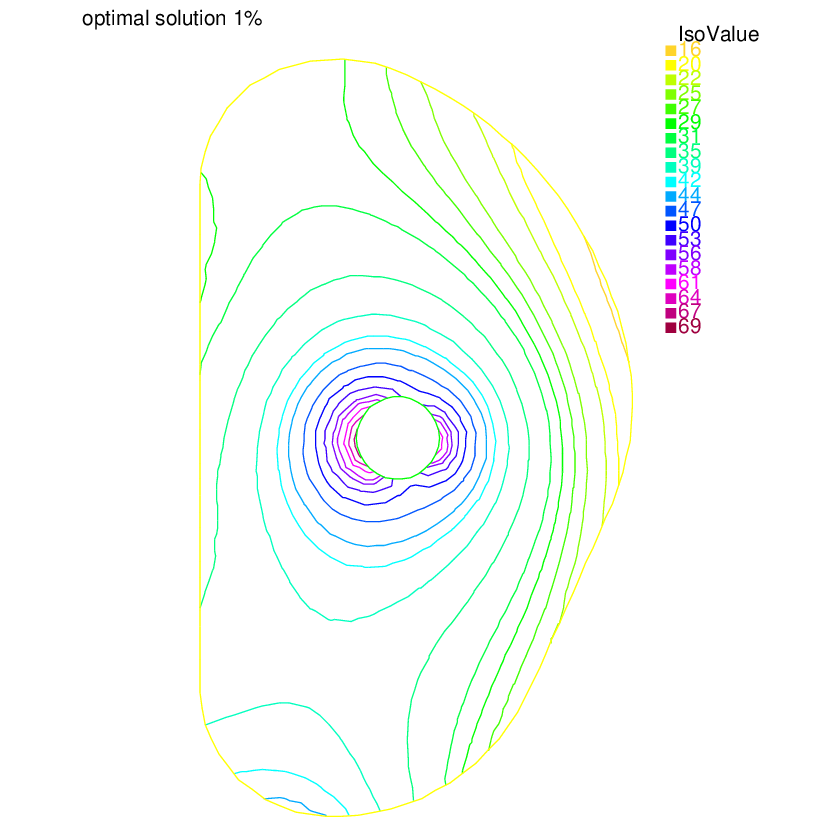}&
\includegraphics[height=6cm,angle=0]{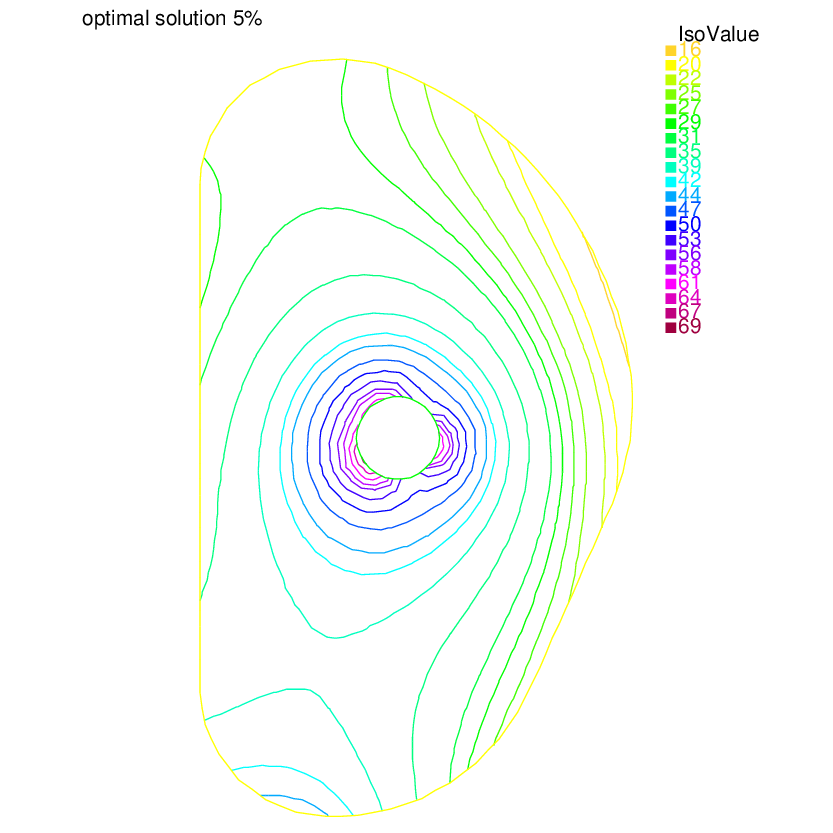}\\
\end{tabular}
\caption{\label{fig:twinuvar} First test case (TC1), $u_{ref}$ given by Eq. (\ref{eqn:uref1}). 
Top left: reference solution $\psi_{N,ref}(u_{ref},g) $. Top right: recovered solution with no noise on the data. Bottom left: recovered solution with a $1\%$ noise on the data. Bottom right: recovered solution with a $5\%$ noise.}
\end{center}
\end{figure}

\begin{figure}
\begin{center}
\begin{tabular}{ll}
\includegraphics[height=6cm,angle=0]{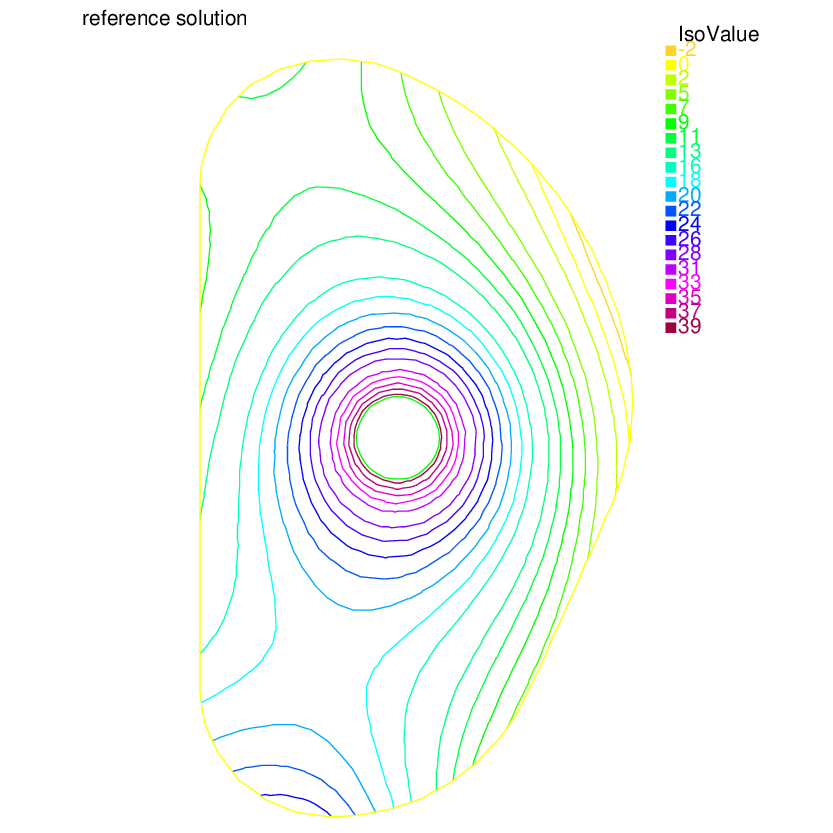}& 
\includegraphics[height=6cm,angle=0]{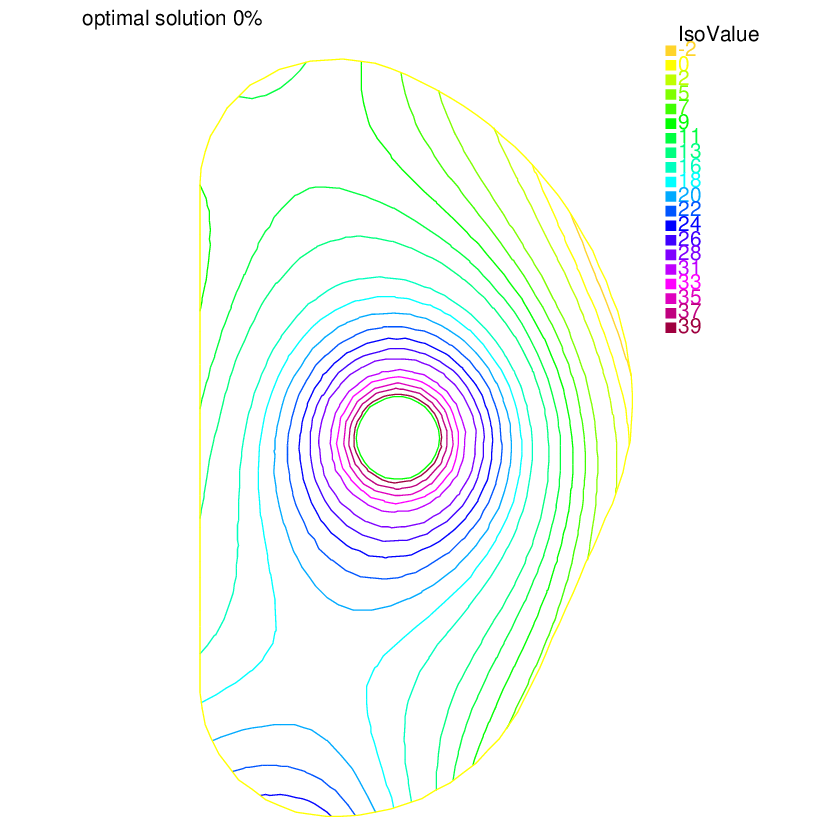}\\ 
\includegraphics[height=6cm,angle=0]{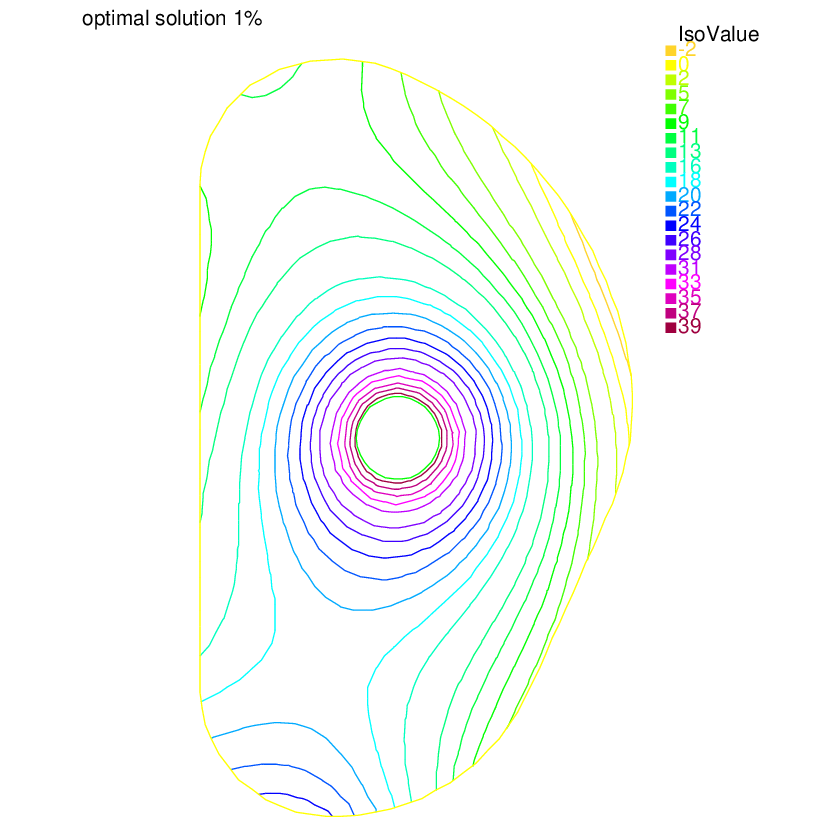}&
\includegraphics[height=6cm,angle=0]{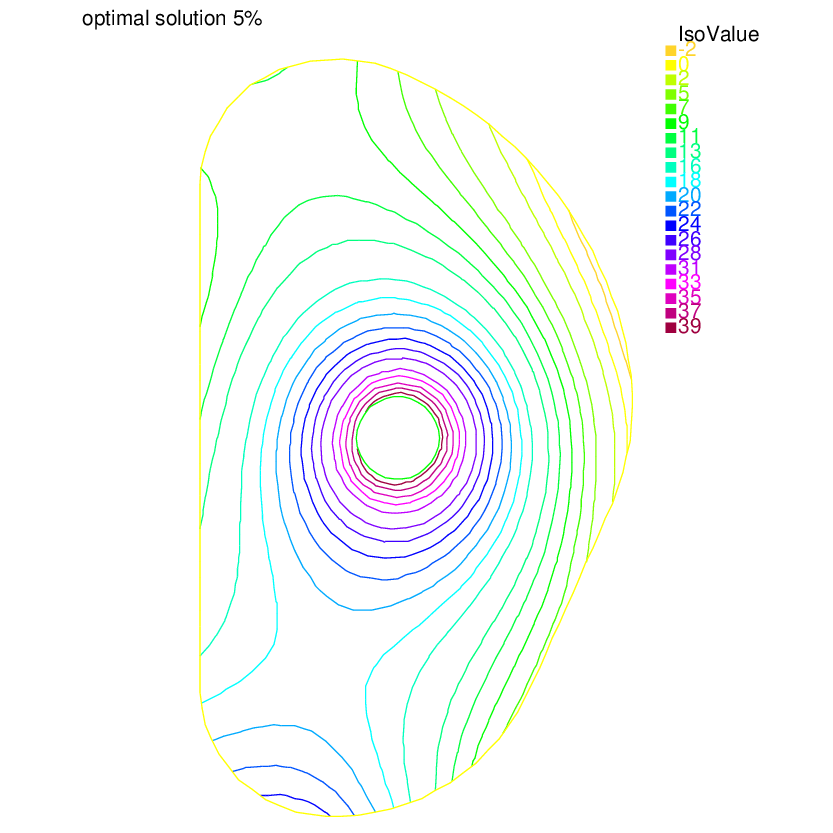}\\
\end{tabular}
\caption{\label{fig:twinucst}
Second test case (TC2), $u_{ref}$ given by Eq. (\ref{eqn:uref2}). 
Top left: reference solution $\psi_{N,ref}(u_{ref},g) $. Top right: recovered solution with no noise on the data. Bottom left: recovered solution with a $1\%$ noise on the data. Bottom right: recovered solution with a $5\%$ noise.
}
\end{center}
\end{figure}

\begin{figure}
\begin{center}
\begin{tabular}{ll}
\includegraphics[height=5cm,angle=0]{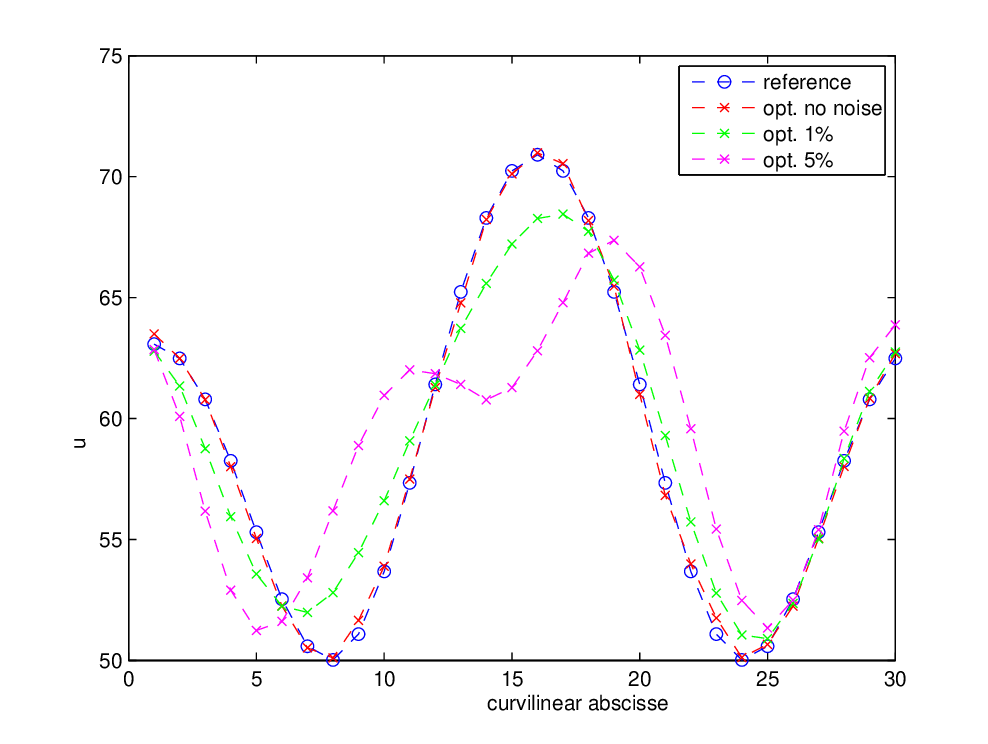}&
\includegraphics[height=5cm,angle=0]{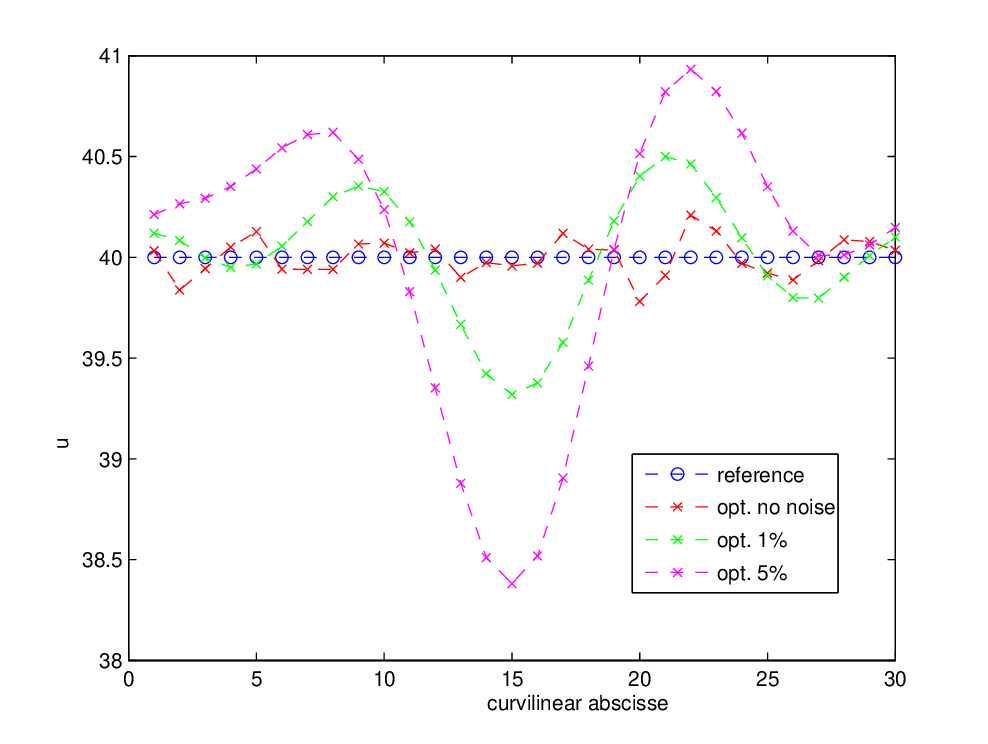}
\end{tabular}
\caption{\label{fig:twinucstuvar}$u_{ref}$ and the recovered $u_{opt}$ 
for the 3 levels of noise on the data. Left: TC1. Right TC2.}
\end{center}
\end{figure}

\begin{figure}
\begin{center}
\includegraphics[height=6cm,angle=0]{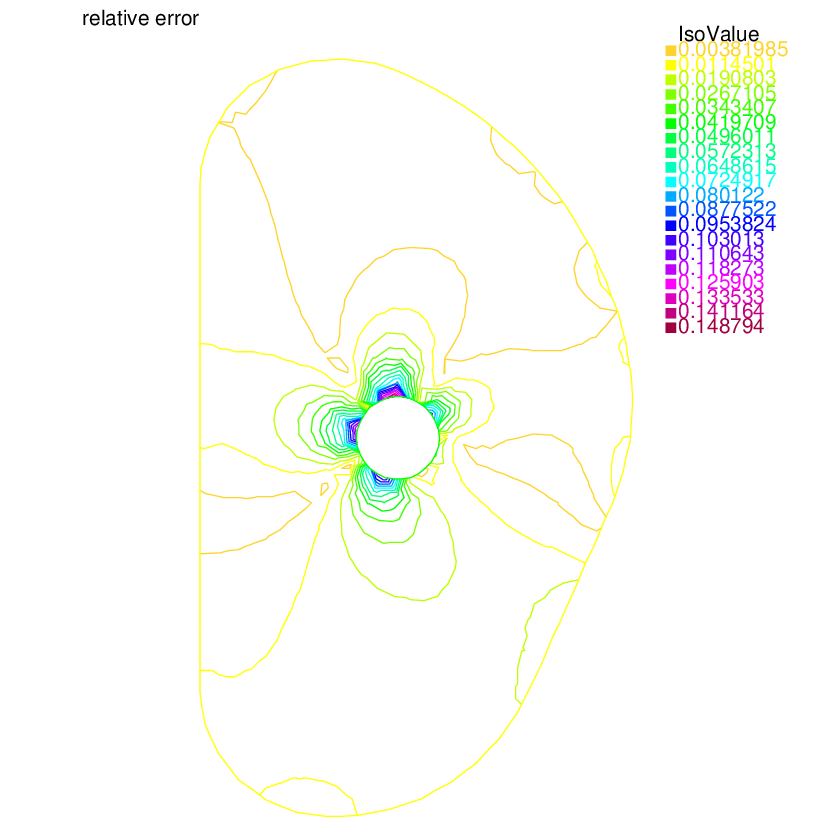}\\
\caption{\label{fig:twinerror} Relative error $|\psi_{opt}-\psi_{opt}|/|\psi_{ref}|$ for TC1 
with $5\%$ noise.}
\end{center}
\end{figure}

The results are shown on Figs. \ref{fig:twinuvar} and \ref{fig:twinucst} 
where the reference and recovered solutions are shown for the three levels of noise considered.
The results are excellent for the noise free case in which the Dirichlet boundary condition $u$ is almost perfecty recovered (Fig. \ref{fig:twinucstuvar}). The differences between $u_{opt}$ and $u_{ref}$ increase with the level of noise (Fig. \ref{fig:twinucstuvar} and Tab. \ref{tab:twinerroru}).
As it is often the case in this type of inverse problems 
the most important errors on $\psi_{opt}$ are localized close to the boundary 
$\Gamma_I$ and vanishes as we move away from it (Fig. \ref{fig:twinerror}).

Tables \ref{tab:twinuvar} and \ref{tab:twinucst} sumarize the evolution of the values of $J$, $R_D$ and $J_\varepsilon$ for the different noise level. First guess values ($u=0$) are also provided for comparison. 
Please note that the regularization parameter 
was chosen differently from one experiment to another depending on the noise level. This was tuned by hand. In the next section we propose to use the L-curve method \cite{Hansen:1998} 
to choose the value of $\varepsilon$.

\begin{table}
\centering
\begin{tabular}{|c|c|c|}
\hline\\[-10pt]
     noise level                 & error  TC1 & error TC2 \\
\hline\\[-10pt]
$0\%$           &  $0.0131$   &  $0.0055$ \\
$1\%$              &  $0.0659$   &  $0.0170$ \\
$5\%$              &  $0.1526$   &  $0.0405$ \\[1pt]
\hline
\end{tabular}
\caption{\label{tab:twinerroru}Maximum relative error $\ds \frac{|u_{opt}-u_{ref}|}{|u_{ref}|}$ 
for TC 1 and 2 }. 
\end{table}

\begin{table}[h!]
\centering
\begin{tabular}{|c|cccc|}
\hline\\[-10pt]
                      & $J$ & $R_D$ & $J_\varepsilon$ & $\varepsilon$ \\
\hline\\[-10pt]
$u=0\ \mathrm{no\ noise}$                 &  $46.8643$   &    $0$   &   $46.8643$          &            \\
$u_{opt}\ \mathrm{no\ noise} $ & $0.0021$  & $46.8722$  &  $0.0026$ & $10^{-5}$  \\
$u_{opt}\ 1\%\ \mathrm{noise} $ & $1.8443$ & $46.5553$ & $1.8676$ & $5 \times 10^{-4}$	\\
$u_{opt}\ 5\%\ \mathrm{noise} $ &  $9.2180$   &  $46.5575$ &  $9.2646$     &    $10^{-3}$ \\[1pt]
\hline
\end{tabular}
\caption{\label{tab:twinuvar} TC1 results. Values of the error functional, the regularization term, the total cost function and the chosen regularization parameter for 
the initial guess (row 1), the optimal solutions for different noise levels (row 2, 3 and 4). }
\end{table}

\begin{table}[h!]
\centering
\begin{tabular}{|c|cccc|}
\hline\\[-10pt]
                      & $J$ & $R_D$ & $J_\varepsilon$ & $\varepsilon$ \\
\hline\\[-10pt]
$u=0\ \mathrm{no\ noise}$         &  $30.7231$   &    $0$   &   $30.7231$   &    \\
$u_{opt}\ \mathrm{no\ noise} $ & $0.0003$  & $30.7242$  &  $0.0006$ & $10^{-5}$  \\
$u_{opt}\ 1\%\ \mathrm{noise} $ & $0.7300$ & $30.7159$ & $0.7607$ & $10^{-3}$	\\
$u_{opt}\ 5\%\ \mathrm{noise} $ &  $3.6516$   &  $30.6822$ &  $3.8050$ & $5\times10^{-3}$ \\[1pt]
\hline
\end{tabular}
\caption{\label{tab:twinucst} TC2 results. 
Values of the error functional, the regularization term, 
the total cost function and the chosen regularization parameter for 
the initial guess (row 1), the optimal solutions for different noise levels (row 2, 3 and 4). }
\end{table}

\clearpage

\subsection{An ITER equilibrium}

\begin{figure}
\begin{center}
\includegraphics[height=5cm,angle=0]{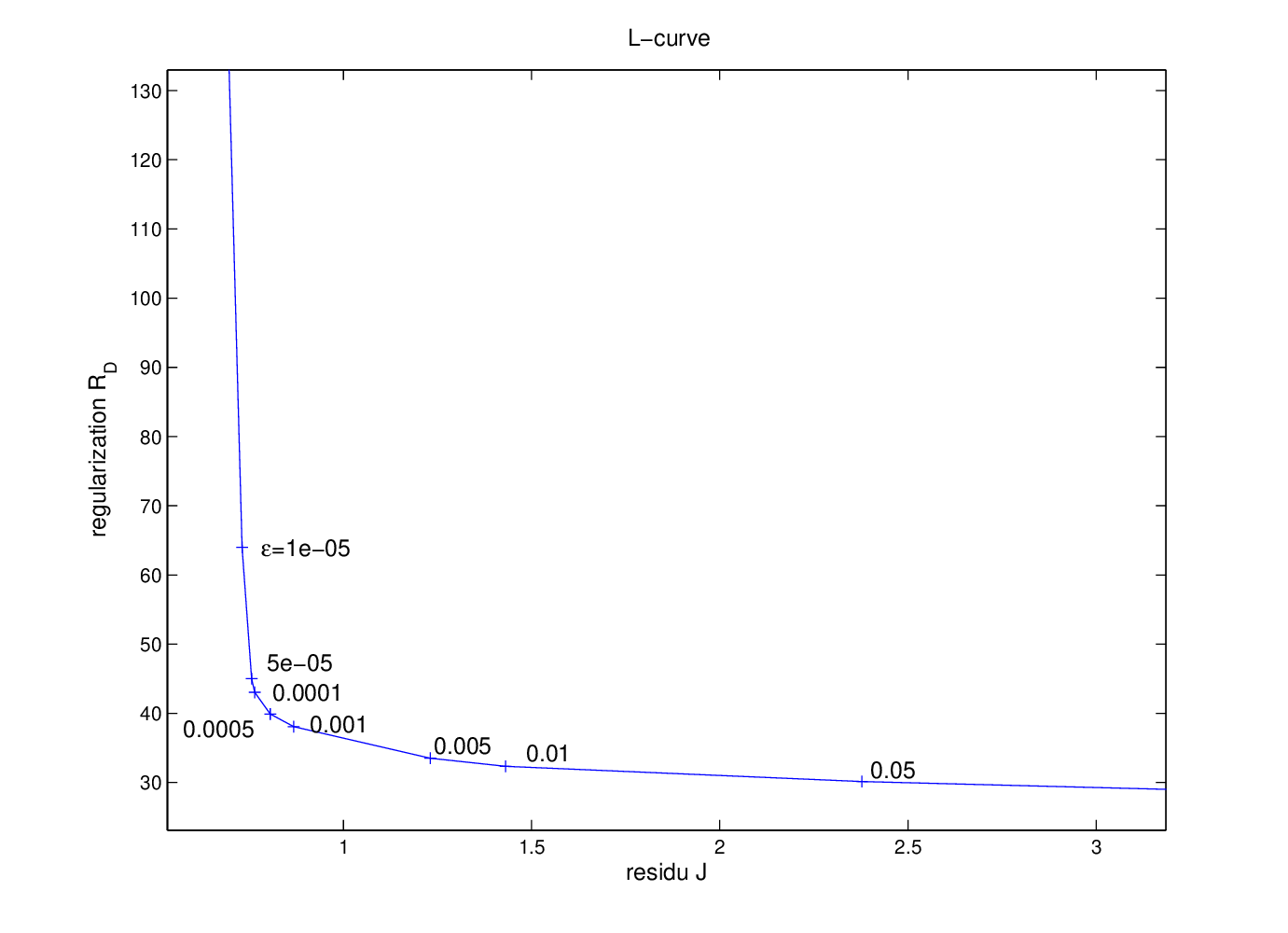}
\end{center}
\caption{\label{fig:lcurve} L-curve computed for the ITER case. The corner is located at $\varepsilon=5\times10^{-4}$.}
\end{figure}

In this last numerical experiment we consider a 'real' ITER case. 
Measurements of the magnetic field are provided by the plasma 
equilibrium code CEDRES++ \cite{GrandGirard:1999}. 
These mesurements are interpolated to provide $f$ and $g$ on $\Gamma_V$. 
The regularized error functional is then minimized to compute the optimal $u_{opt}$. 
The choice of the regularization parameter $\varepsilon$ is made thanks to the computation of the 
L-curve shown on Fig. \ref{fig:lcurve}. It is a plot of $(J(u_{opt})(\varepsilon),R_D(u_{opt})(\varepsilon))$ 
as $\varepsilon$ varies. 
The corner of the L-shaped curve provides a value of $\varepsilon=5.10^{-4}$.

\begin{figure}
\begin{center}
\includegraphics[height=5cm,angle=0]{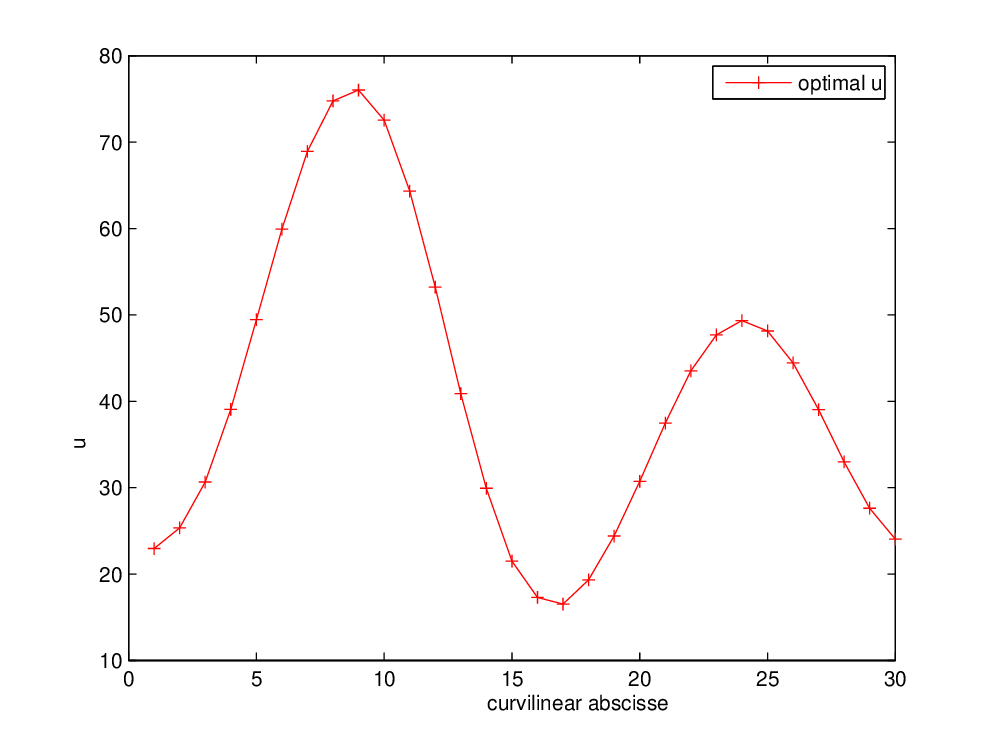}
\end{center}
\caption{\label{fig:uopt}Optimal $u_{opt}$ for the ITER case.}
\end{figure}

\begin{figure}
\begin{center}
\includegraphics[height=7cm,angle=0]{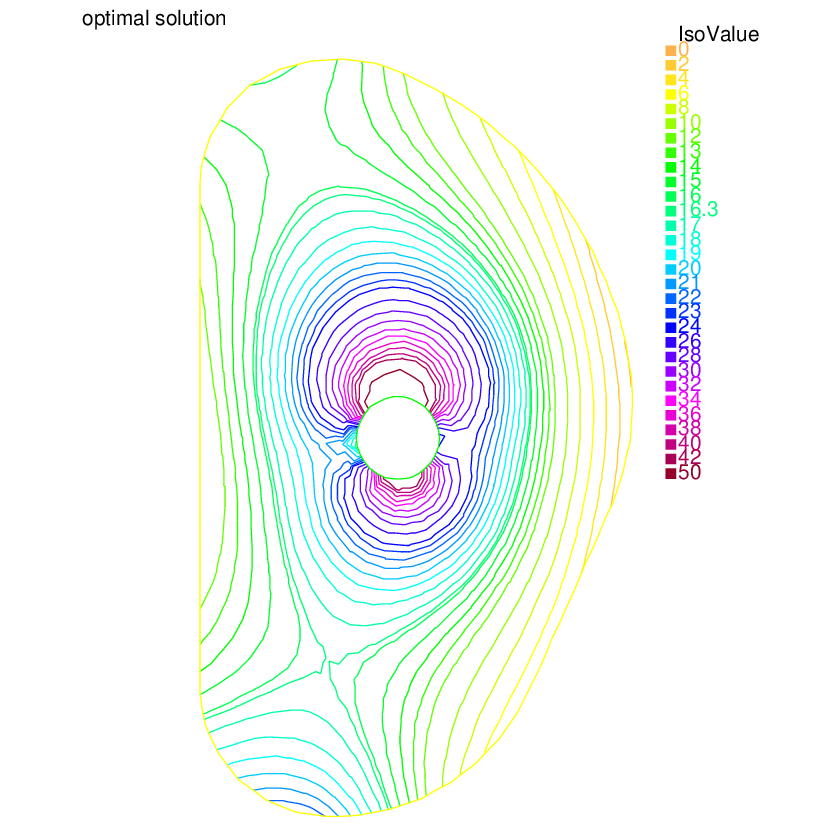}
\end{center}
\caption{\label{fig:optimalsolutioniter}Optimal solution for the ITER case. 
The plasma is in an X-point configuration}
\end{figure}

The computed $u_{opt}$ is shown on Fig. \ref{fig:uopt} and numerical 
values are given in Tab. \ref{tab:iter}. 
The recovered poloidal flux $\psi$ is shown on Fig. \ref{fig:optimalsolutioniter}. 
The boundary of the plasma is found to be the isoflux $\psi=16.3$ which 
shows an X-point configuration.

\begin{table}[h!]
\centering
\begin{tabular}{|c|cccc|}
\hline\\[-10pt]
                      & $J$ & $R_D$ & $J_\varepsilon$ & $\varepsilon$ \\
\hline\\[-10pt]
$u=0$         &  $31.1026$   &  $0$ &   $31.1026$   &    \\
$u_{opt}$ & $0.8053$  & $39.9169$  &  $0.8253$ & $5\times10^{-4}$ \\
\hline
\end{tabular}
\caption{\label{tab:iter}ITER case results. Values of the error functional, 
the regularization term, the total cost function and the chosen regularization parameter for 
the initial guess (row 1) and the optimal solution (row 2)}
\end{table}

\clearpage

\section{Conclusion}

We have presented a numerical method for the computation of the poloidal flux in the vacuum region 
surrounding the plasma in a Tokamak. The algorithm is based on the optimization of a regularized error functional.
This computation enables in a second step the identification of the plasma boundary.

Numerical experiments have been conducted. 
They show that the method is precise and robust to noise on the Cauchy input data. 
It is fast since the optimization reduces to the resolution of a linear system 
of very reasonable dimension. Successive equilibrium reconstructions can be conducted 
very rapidly since the matrix of this linear system can be completely precomputed 
and only the right hand side has to be updated. 
The L-curve method proved to be efficient to specify the regularization parameter.

\clearpage

\section*{Appendix. Proof of Proposition 1}

{\it 1.} We need to prove the continuity and the coercivity of $J_{\varepsilon}$.\\
Continuity.\\
The maps $v \mapsto \psi_D(v)$ and $v \mapsto \psi_N(v)$ are 
continuous and linear from 
$H^{1/2}(\Gamma_I)$ to $H^1(\Omega)$. Moreover since $\tilde{\psi}_D(f)$ and $\tilde{\psi}_N(g)$ are in 
$H^1(\Omega)$ 
and $ r_M \ge r \ge r_m>0$ in $\Omega$ it is shown with Cauchy Schwarz  that the bilinear forms $s_D$ and $s_N$, the linear form $l$ and thus $J_\varepsilon$ are continuous on $H^{1/2}(\Gamma_I)$.\\

Coercivity.\\
The bilinear form $s_D$ is coercive on $H^{1/2}(\Gamma_I)$. One obtains this from the fact that 
$\psi_D(v) \in H^1_0(\Omega,\Gamma_V)$ 
and the Poincar\'e inequality holds, and from the continuity of the application 
$\psi_D(v) \in H^1(\Omega) \rightarrow \psi_D(v)|_{\Gamma_I}=v \in H^{1/2}(\Gamma_I)$. 

On the contrary, since for $\psi_N(v) \in H^1(\Omega)$ the seminorm does not bound the $L^2$ norm, 
the bilinear form $s_N$ is not coercive and because of the minus sign in $s=s_D - s_N$ 
we need to prove that $s(v,v)\ge0$ to obtain the coercivity 
of the bilinear part of functional $J_\varepsilon$. One can use the same type of argument 
as in \cite{BenBelgacem:2005} to de so.

Eventually it holds that 
$$
\ds \frac{1}{2}s(v,v)+ \frac{\varepsilon}{2}s_D(v,v) \ge C \varepsilon || v ||^2_{H^{1/2}(\Gamma_I)}
$$   
Using the continuity and the coercivity of $J_\varepsilon$ it results from \cite{Lions2:1968} 
that problem $P_\varepsilon$ admits a unique solution $u_\varepsilon \in \mathcal{U}$.

The solution $u_\varepsilon$ is characterized by the first order optimality condition 
which is written as the following well-posed variational problem
\begin{equation}
\label{eqn:opt}
(J'_\varepsilon(u_\varepsilon),v)=\varepsilon s_D(u_\varepsilon,v) + s_D(u_\varepsilon,v) - s_N(u_\varepsilon,v) 
-l(v)=0\quad \forall v \in \mathcal{U}
\end{equation}

which as in Eq. (\ref{eqn:normal-boundary}) can be understood as an 
equality on $\Gamma_I$.

~\\
{\it 2.} The stability result is deduced from the optimality condition (\ref{eqn:opt}). 

Let $u^1_\varepsilon$ (resp. $u^2_\varepsilon$) be the solution associated to $(f^1,g^1)$ (resp. $(f^2,g^2)$). 
Substracting the two optimality conditions, 
choosing $v=u^1_\varepsilon-u^2_\varepsilon$ and using the coercivity leads to
$$
C \varepsilon || u^1_\varepsilon-u^2_\varepsilon||^2_{H^{1/2}(\Gamma_I)} \le | (l_1-l_2)(u^1_\varepsilon-u^2_\varepsilon)|
$$
The map $f \mapsto \tilde{\psi}_D(f)$ is linear and continuous 
from $H^{1/2}(\Gamma_V)$ to $H^1(\Omega)$, 
and so is the map  $g \mapsto \tilde{\psi}_N(g)$ 
from $H^{-1/2}(\Gamma_V)$ to $H^1(\Omega)$. 
Using these facts and Cauchy Schwarz it follows that
$$
 || u^1_\varepsilon-u^2_\varepsilon||_{H^{1/2}(\Gamma_I)} \le 
 \ds \frac{C'}{r_m C}\frac{1}{\varepsilon} (||f^1-f^2 ||_{H^{1/2}(\Gamma_V)} + ||g^1-g^2 ||_{H^{-1/2}(\Gamma_V)})
$$

{\it 3.} For this point the proof of Proposition $3.2$ in \cite{Azaiez:2006} can be adpated. 
A sketch of the proof is as follows. Let us suppose that there exists $u \in \mathcal{U}$ 
such that $\psi_D(u,f)=\psi_N(u,g)$. A key point is to show that 
$s_D(u_\varepsilon,u_\varepsilon) \rightarrow s_D(u,u)$ when $\varepsilon \rightarrow 0$. 
Then in a second step using the optimality conditions 
for $u$ and $u_\varepsilon$ it is shown that
$$
s_D( u_\varepsilon-u,  u_\varepsilon-u) 
\le s_D(u, u) - s_D( u_\varepsilon,u_\varepsilon)
$$ 
which gives the result thanks to the coercivity of $s_D$ in $H^{1/2}(\Gamma_I)$. 

\clearpage

\bibliographystyle{plain}
\bibliography{biblio-vacuumflux}

\begin{thebibliography}{10}

\bibitem{Andrieux:2006}
S.~Andrieux, T.N. Baranger, and A.~Ben~Abda.
\newblock Solving cauchy problems by minimizing an energy-like functional.
\newblock {\em Inverse Problems}, 22:115--133, 2006.

\bibitem{Andrieux:2005}
S.~Andrieux, A.~Ben~Abda, and T.N. Baranger.
\newblock Data completion via an energy error functional.
\newblock {\em C.R. Mecanique}, 333:171--177, 2005.

\bibitem{Azaiez:2006}
M.~Azaiez, F.~Ben~Belgacem, and H.~El~Fekih.
\newblock On {C}auchy's problem: {II}. {C}ompletion, regularization and
  approximation.
\newblock {\em Inverse Problems}, 22:1307--1336, 2006.

\bibitem{BenAbda:2009}
A~Ben~Abda, M.~Hassine, M.~Jaoua, and M.~Masmoudi.
\newblock Topological sensitivity analysis for the location of small cavities
  in stokes flows.
\newblock {\em SIAM J. Cont. Opt.}, 2009.

\bibitem{BenBelgacem:2005}
F.~Ben~Belgacem and H.~El~Fekih.
\newblock On {C}auchy's problem: I. {A} variational {S}teklov-{P}oincar\'e
  theory.
\newblock {\em Inverse Problems}, 21:1915--1936, 2005.

\bibitem{Blum:1989}
J.~Blum.
\newblock {\em Numerical Simulation and Optimal Control in Plasma Physics with
  Applications to Tokamaks}.
\newblock Series in Modern Applied Mathematics. Wiley Gauthier-Villars, Paris,
  1989.

\bibitem{Bonnet:2005}
M.~Bonnet and A.~Constantinescu.
\newblock Inverse problems in elasticity.
\newblock {\em Inverse Problems}, 21(2), 2005.

\bibitem{Bourgeois:2010}
L.~Bourgeois and J.~Dard\'e.
\newblock A quasi-reversibility approach to solve the inverse obstacle problem.
\newblock {\em Inverse Problems and Imaging}, 4/3:351--377, 2010.

\bibitem{Braams:1991}
B.J. Braams.
\newblock The interpretation of tokamak magnetic diagnostics.
\newblock {\em Nuc. Fus.}, 33(7):715--748, 1991.

\bibitem{Chaabane:1999}
S.~Chaabane and M.~Jaoua.
\newblock Identification of robin coefficients by means of boundary
  measurements.
\newblock {\em Inverse Problems}, 15(6):1425, 1999.

\bibitem{Ciarlet:1980}
P.G. Ciarlet.
\newblock {\em The Finite Element Method For Elliptic Problems}.
\newblock North-Holland, 1980.

\bibitem{Courant:1962}
R.~Courant and D.~Hilbert.
\newblock {\em Methods of Mathematical Physics}, volume 1-2.
\newblock Interscience, 1962.

\bibitem{Feneberg:1984}
W.~Feneberg, K.~Lackner, and P.~Martin.
\newblock Fast control of the plasma surface.
\newblock {\em Computer Physics Communications}, 31(2):143--148, 1984.

\bibitem{Fischer:2011}
Y.~Fischer.
\newblock {\em Approximation dans des classes de fonctions analytiques
  g\'en\'eralis\'ees et r\'esolution de probl\`emes inverses pour les
  tokamaks}.
\newblock Phd thesis, Universit\'e de Nice Sophia Antipolis, France, 2011.

\bibitem{Fischer:2012}
Y.~Fischer, B.~Marteau, and Y.~Privat.
\newblock Some inverse problems around the tokamak tore supra.
\newblock {\em Comm. Pure and Applied Analysis}, to appear.

\bibitem{Grad:1958}
H.~Grad and H.~Rubin.
\newblock Hydromagnetic equilibria and force-free fields.
\newblock In {\em 2nd U.N. Conference on the Peaceful uses of Atomic Energy},
  volume~31, pages 190--197, Geneva, 1958.

\bibitem{GrandGirard:1999}
V.~Grandgirard.
\newblock {\em Mod\'elisation de l'\'equilibre d'un plasma de tokamak - Tokamak
  plasma equilibrium modelling}.
\newblock Phd thesis, Universit\'e de Besan\c{c}on, France, 1999.

\bibitem{Haddar:2005}
H.~Haddar and R.~Kress.
\newblock Conformal mappings and inverse boundary value problem.
\newblock {\em Inverse Problems}, 21:935--953, 2005.

\bibitem{Hansen:1998}
C.~Hansen.
\newblock {\em Rank-Deficient and Discrete Ill-Posed Problems: Numerical
  Aspects of Linear Inversion}.
\newblock SIAM, Philadelphia, 1998.

\bibitem{Kohn:1990}
R.V. Kohn and A.~McKenney.
\newblock Numerical implementation of a variational method for electrical
  impedance tomography.
\newblock {\em Inverse Problems}, 6(3):389, 1990.

\bibitem{Kohn:1985}
R.V. Kohn and M.S. Vogelius.
\newblock Determining conductivity by boundary measurements: {II}. {I}nterior
  results.
\newblock {\em Commun. Pure Appl. Math.}, 31:643--667, 1985.

\bibitem{Kohn:1989}
R.V. Kohn and M.S. Vogelius.
\newblock Relaxation of a variational method for impedance computed tomography.
\newblock {\em Commun. Pure Appl. Math.}, 11:745--777, 1987.

\bibitem{Lackner:1976}
K.~Lackner.
\newblock Computation of ideal {MHD} equilibria.
\newblock {\em Computer Physics Communications}, 12:33--44, 1976.

\bibitem{Ladeveze:1983}
P.~Ladeveze and D.~Leguillon.
\newblock Error estimate procedure in the finite element method and
  applications.
\newblock {\em SIAM J. Num. Anal.}, 20(3):485--509, 1983.

\bibitem{Lions2:1968}
J.L. Lions.
\newblock {\em Contr\^ole optimal de syst\`emes gouvern\'es par des \'equations
  aux d\'eriv\'ees partielles (Optimal control of systems governed by partial
  differential equations)}.
\newblock Dunod, Paris, 1968.

\bibitem{OBrien:1993}
D.P. O'Brien, J.J Ellis, and J.~Lingertat.
\newblock Local expansion method for fast plasma boundary identification in
  {JET}.
\newblock {\em Nuc. Fus.}, 33(3):467--474, 1993.

\bibitem{Quarteroni:1999}
A.~Quarteroni and V~Alberto.
\newblock {\em Domain Decomposition Methods for Partial Differential
  Equations}.
\newblock Oxford University Press, 1999.

\bibitem{Saint-Laurent:2001}
F.~Saint-Laurent and G.~Martin.
\newblock Real time determination and control of the plasma localisation and
  internal inductance in {T}ore {S}upra.
\newblock {\em Fusion Engineering and Design}, 56-57:761--765, 2001.

\bibitem{Sartori:2003}
F.~Sartori, A.~Cenedese, and F.~Milani.
\newblock {JET} real-time object-oriented code for plasma boundary
  reconstruction.
\newblock {\em Fus. Engin. Des.}, 66-68:735--739, 2003.

\bibitem{Shafranov:1958}
V.D. Shafranov.
\newblock On magnetohydrodynamical equilibrium configurations.
\newblock {\em Soviet Physics JETP}, 6(3):1013, 1958.

\bibitem{Shafranov:1972}
V.D. Shafranov and L.E. Zakharov.
\newblock Use of the virtual-casing principle in calculating the containing
  magnetic field in toroidal plasma systems.
\newblock {\em Nuc. Fus.}, 12:599--601, 1972.

\bibitem{Wesson:2004}
J.~Wesson.
\newblock {\em Tokamaks}, volume 118 of {\em International Series of Monographs
  on Physics}.
\newblock {O}xford {U}niversity {P}ress {I}nc., {N}ew {Y}ork, Third Edition,
  2004.

\end{thebibliography}

\end{document}